\documentclass[11pt]{article}
\usepackage{latexsym}
\usepackage{amsfonts}
\usepackage{amssymb}

\newcommand{\be}{\begin{equation}}
\newcommand{\ee}{\end{equation}}
\newcommand{\bea}{\begin{eqnarray}}
\newcommand{\eea}{\end{eqnarray}}
\newcommand{\bean}{\begin{eqnarray*}}
\newcommand{\eean}{\end{eqnarray*}}
\newcommand{\brray}{\begin{array}}
\newcommand{\erray}{\end{array}}

\newcommand{\newsection}[1]{\setcounter{equation}{0} \setcounter{dfn}{0}
\section{#1}}

\newtheorem{dfn}{Definition}[section]
\newtheorem{thm}[dfn]{Theorem}
\newtheorem{lmma}[dfn]{Lemma}
\newtheorem{ppsn}[dfn]{Proposition}
\newtheorem{crlre}[dfn]{Corollary}
\newtheorem{xmpl}[dfn]{Example}
\newtheorem{rmrk}[dfn]{Remark}

\newcommand{\bdfn}{\begin{dfn}}
\newcommand{\bthm}{\begin{thm}}
\newcommand{\blmma}{\begin{lmma}}
\newcommand{\bppsn}{\begin{ppsn}}
\newcommand{\bcrlre}{\begin{crlre}}
\newcommand{\bxmpl}{\begin{xmpl}}
\newcommand{\brmrk}{\begin{rmrk}}

\newcommand{\edfn}{\end{dfn}}
\newcommand{\ethm}{\end{thm}}
\newcommand{\elmma}{\end{lmma}}
\newcommand{\eppsn}{\end{ppsn}}
\newcommand{\ecrlre}{\end{crlre}}
\newcommand{\exmpl}{\end{xmpl}}
\newcommand{\ermrk}{\end{rmrk}}

\newcommand{\IC}{\mathbb{C}}
\newcommand{\IZ}{\mathbb{Z}}

\newcommand{\afa}{\alpha}
\newcommand{\bta}{\beta}
\newcommand{\eps}{\epsilon}

\newcommand{\cla}{{\cal A}}
\newcommand{\clb}{{\cal B}}
\newcommand{\clh}{{\cal H}}
\newcommand{\clk}{{\cal K}}

\def \bbt {\mbox{\boldmath $t$}}

\newcommand{\prf}{\noindent{\it Proof\/}: }

\newcommand{\raro}{\rightarrow}

\newcommand{\seq}{\subseteq}

\newcommand{\lgl}{\langle}
\newcommand{\rgl}{\rangle}

\newcommand{\NI}{\noindent}
\newcommand {\CC}{\centerline}
\def \qed { \mbox{}\hfill $\Box$\vspace{1ex}}

\begin{document}

\author{{\large Arupkumar Pal\thanks{Partially supported by the 
    Jawaharlal Nehru Centre for Advanced Scientific Research, 
    Bangalore, India.}}\\
         Indian Statistical Institute,\\[-.5ex]
         7, SJSS Marg, New Delhi--110\,016, INDIA\\[-.5ex]
         email: arup@isid.ac.in}
\title{Regularity of  Operators on \\Essential Extensions of the Compacts}
\date{May\,25, 1998\\ Revised: October\,22, 1998}   
\maketitle
   \begin{quotation}
\CC{\bf Abstract}

A semiregular operator on a Hilbert $C^*$-module, or equivalently,
on the $C^*$-algebra of `compact' operators on it, is a closable
densely defined operator whose adjoint is also densely defined. 
It is shown that for operators on  extensions of compacts by unital or abelian
$C^*$-algebras, semiregularity leads to regularity. Two examples 
coming from quantum groups are discussed.\\[2ex]
{\bf AMS Subject Classification No.:} {\large 46}H{\large 25},
                                  {\large 47}C{\large 15}
   \end{quotation}

\newsection{Introduction}
Hilbert $C^*$-modules arise in many different areas, for 
example, in the study of locally compact quantum groups and
their representations, in KK-theory, in noncommutative geometry,
and in the study of completely positive maps between $C^*$-algebras.
A regular operator on a Hilbert $C^*$-module is an analogue of 
a closed operator on a Hilbert space that naturally arises
in many of the above contexts. A closed and densely defined operator $T$
on a Hilbert $C^*$-module $E$ is called regular if its adjoint $T^*$
is also densely defined and if the range of $(I+T^*T)$ is dense in $E$.
Every regular operator on a Hilbert $C^*$-module $E$ is uniquely 
determined by a (bounded) adjointable operator on $E$, called its
$z$-transform. This fact is exploited when dealing with regular operators,
as the adjointable operators, being bounded, are more easily 
manageable than unbounded operators. But given an unbounded
operator, the first and the most basic problem is to decide whether
or not it is regular. In \cite{wo1}, Woronowicz investigated this 
problem using graphs of operators, and proved a few results
(see proposition~2.2, theorem~2.3 and examples~1--3 in \cite{wo1}).
In particular, he was able to conclude the regularity of some very 
simple functions of a regular operator $T$, 
like $T+a$ where $a$ is an adjointable operator, and  $Ta$ and $aT$ 
where $a$ is an invertible adjointable operator.

  The problem was later attacked from a different angle in \cite{pa2}.
A somewhat larger class of operators, called the {\em semiregular}
operators were considered. A semiregular operator is a closable
densely defined operator whose adjoint is also densely defined.
Though regularity is quite difficult to ascertain, semiregularity
is not. The problem then investigated in \cite{pa2} was
`when is a semiregular operator regular?'. The first step was to 
reduce the problem to a problem on $C^*$-algebras by establishing 
that semiregular operators on a Hilbert $C^*$-module $E$ correspond,
in a canonical manner, to those on the $C^*$-algebra $\clk(E)$ of
`compact' operators on $E$. The question to be answered next is then
`for what class of $C^*$-algebras is a closed semiregular operator
regular (or admits regular extension)?' Among other things, it was 
established that for abelian $C^*$-algebras as well as for subalgebras 
of $\clb_0(\clh)$, closed semiregular operators are indeed regular. 
In the present paper, we will extend the results to a class of $C^*$-algebras
that contain $\clb_0(\clh)$ as an essential ideal. Most of the 
results, however, are valid in a more general situation where 
$\clb_0(\clh)$ is 
replaced by any essential ideal $K$. Since it involves almost no 
extra work, the results are stated in this general set up.
In section~2, we develop the necessary background for proving the main 
results which are presented in section~3. Finally in section~4, we discuss 
two examples that arise in the context of quantum groups and are covered 
by the results in section~3. We have assumed elements of $C^*$-algebra
theory and Hilbert $C^*$-module theory as can be found, for example,
in Pedersen (\cite{pe}) and Lance (\cite{la}) respectively.

Now, why are essential extensions of the compacts important 
in the context of the problem? Firstly, because they cover examples
that arise naturally, like the quantum complex plane which has been 
discussed later in this paper. Secondly and perhaps more importantly,
because they arise as irreducible representations of all type~I
$C^*$-algebras. For a large class of type~I $C^*$-algebras, one would 
be able to conclude by the results here that all irreducible `fibres'
of a semiregular operator $S$ are regular. This fact, along with some mild
restrictions on $S$ should then lead to its regularity.

\NI{\bf Notations.} We will follow standard notations mostly. 
So, for example, $\clh$ is
a complex separable Hilbert space, $\clb_0(\clh)$ is the algebra of
compact operators on $\clh$; $\cla$ is a $C^*$-algebra, $M(\cla)$
and $LM(\cla)$ are the space of multipliers and left multipliers
respectively of $\cla$. For a topological space $X$, $C_0(X)$ will
denote the $C^*$-algebra of continuous functions on $X$ vanishing at
infinity. The $C^*$-algebra $\cla$ that we will primarily
be interested in, {\it will always be assumed to be separable} (this of course
will not be true for all $C^*$-algebras that we deal with; for example 
the multiplier algebra of a nonunital $C^*$-algebra is 
never separable).

\newsection{Restriction to an Ideal}

Let $\cla$ be a nonunital $C^*$-algebra and let $K$ be an essential ideal
in $\cla$. Since $\cla$ is essential in $M(\cla)$, it follows that $K$ is
essential in $M(\cla)$. It is easy to see then that there is an
injective homomorphism from $M(\cla)$ to $M(K)$ through which $M(\cla)$
can be thaught of as sitting inside $M(K)$.

For the rest of this paper, we will always assume that
$K\seq\cla\seq M(\cla)\seq M(K)$.

Before we proceed further, let us recall the definition of a
semiregular operator.

\bdfn[\cite{pa2}] \label{t1.0}
Let $E$ and $F$ be Hilbert $\cla$-modules. An operator $T:E\raro F$ is called
{\em semiregular} if \\[-5ex]
\begin{enumerate}\renewcommand{\theenumi}{\alph{enumi}}
\item $D_T$ is a dense submodule of $E$ {\rm (i.e.\ $D_T\cla\seq D_T$)},
      \vspace{-2ex}
\item $T$ is closable, \vspace{-2ex}
\item $T^*$ is densely defined.
\end{enumerate}
\edfn

Next we list some elementary observations regarding the
restriction of a semiregular operator to an essential ideal.

\bppsn\label{t1.1}
Let $S$ be a closed semiregular operator on $\cla$. Then\\[-5ex]
\begin{enumerate}
\item $D_K:=D(S)\cap K$ is a dense right ideal in $K$,\vspace{-2ex}
\item $S(D_K)\seq K$,\vspace{-2ex}
\item $S_0:=S|_K$ is closed and semiregular,\vspace{-2ex}
\item $D(S)K$ is a core for $S_0$, \vspace{-2ex}
\item $(S|_K)^*=S^*|_K$.
\end{enumerate}
\eppsn

\NI {\it Proof\/}: 
1. That $D_K$ is a right ideal
is obvious. Let us show that it is dense. Choose any $a\in K$. Let
$\{e_\afa\}$ be an approximate identity in $K$. For any $\eps>0$, there 
is an $\tilde{a}\in D(S)$ such that $\|a-\tilde{a}\|<\eps$. Hence
for large enough $\afa$,
\bean
\|\tilde{a}e_\afa-a\| &\leq & \|\tilde{a}-a\|\|e_\afa\|+\|ae_\afa-a\| \\
   &\leq & 2\eps.
\eean
Since $\tilde{a}\in D(S)$, $e_\afa\in K$, $\tilde{a}e_\afa\in D(S)\cap K$.

2. Take an $a\in D_K$. For any $b\in D(S^*)$, $b^*Sa=(S^*b)^*a\in K$. 
Since $D(S^*)$ is dense in $\cla$, we have $b^*(Sa)\in K$ for all $b\in\cla$.
Put $b=Sa$ to get $(Sa)^*(Sa)\in K$. Hence $|Sa|^{1/2}\in K$. 
Now in $\cla$, there exists an element $u$ such that $Sa=u|Sa|^{1/2}$.
Hence $Sa\in K$.

3. For $a\in D(S^*|_K)$ and $b\in D(S|_K)$, we have 
\[
\lgl S^*|_K a,b\rgl =\lgl S^*a,b\rgl=\lgl a,Sb\rgl=\lgl a,S|_Kb\rgl.
\]
Therefore $S^*|_K \seq(S|_K)^*$ and $(S|_K)^*$ is densely defined. 
Now suppose $a_n\in D(S|_K)=J_K$,
and $a_n\raro a$, $S|_K a_n\raro b$. Since $S|_Ka_n=Sa_n$ and $S$ is 
closed, we conclude that $a\in D(S)$ and $Sa=b$. But $a\in K$ also.
Hence $a\in J_K$, and $S|_Ka=b$.

4. Take $a\in D(S|_K)$. If $\{e_\afa\}$ is an approximate identity
for $K$, then $ae_\afa\raro a$ and $S|_K(ae_\afa)=(S|_Ka)e_\afa\raro S|_Ka$.
Since $ae_\afa\in D(S)K$, $D(S)K$ is a core for $S|_K$.

5. We have already seen that $S^*|_K\seq (S|_K)^*$. Let us prove the
reverse inclusion here. For any $a\in D((S|_K)^*)$, $b\in D(S)$,
$k\in K$, we have 
\bean
\lgl a, Sb\rgl k &=& \lgl a, S(bk)\rgl\\
 &=& \lgl (S|_K)^*a,bk\rgl \\
&=&\lgl (S|_K)^*a,b\rgl k.
\eean
Hence $\lgl a, Sb\rgl=\lgl (S|_K)^*a,b\rgl$, so that $a\in D(S^*)$.
Thus $D((S|_K)^*)\seq D(S^*)\cap K=D(S^*|_K)$.\qed

\bppsn\label{t1.2}
Let $S$ and $T$ be semiregular operators on $\cla$ such that
$S|_K=T|_K$. Then\\[-6ex]
\begin{enumerate}
\item $S=T$ on $D(S)\cap D(T)$,\vspace{-2ex}
\item $S^*=T^*$, \vspace{-2ex}
\item if $(S|_K)^{**}=S|_K$, then there exists a maximal closed
semiregular operator on $\cla$ whose restriction to $K$ equals $S|_K$.
\end{enumerate}
\eppsn

\NI {\it Proof\/}: 
1. Take $a\in D(S)\cap D(T)$. For any $k\in K$, 
$ak\in D(S|_K)=D(T|_K)$. Hence $(Sa)k=S(ak)=T(ak)=(Ta)k$.
Therefore $Sa=Ta$.

2. Take any $a\in D(S^*)$, $b\in D(T)$. Then for any $k\in K$, 
\bean
\lgl a, Tb\rgl k &=& \lgl a , T(bk)\rgl\\
&=&\lgl a,S(bk)\rgl\\
&=&\lgl S^*a,bk\rgl\\
&=&\lgl S^*a,b\rgl k.
\eean
Hence $\lgl a, Tb\rgl=\lgl S^*a,b\rgl$. Thus $S^*\seq T^*$. 
Similarly $T^*\seq S^*$.

3. $S^{**}$ is the required operator. 
For, if $T$ is any other semiregular operator whose restriction 
to $K$ is $S|_K$, then $T^*=S^*$, thereby implying $S^{**}=T^{**}$, so that
$T\seq S^{**}$. By part~5 of the forgoing proposition,
$S^{**}|_K=(S^*|_K)^*=(S|_K)^{**}=S|_K$.\qed

Part~3 above tells us, in particular, that if $S|_K$ is regular then 
$S^{**}$ is the maximal semiregular operator on $\cla$ whose restriction 
to $K$ is same as that of $S$.

\blmma\label{t1.3}
If $T$ is regular on $\cla$ with $z$-transform $z$, then $T(K)\seq K$,
and $T|_K$ is a regular operator on $K$ with the same $z$-transform $z$.
\elmma

\NI {\it Proof\/}: 
Observe that $z\in M(\cla)\seq M(K)$,
and $(I-z^*z)^{1/2}K$ contains $(I-z^*z)^{1/2}\cla K =D(T)K$ 
which is dense in $K$.
Hence there exists a regular operator $T_0$ on $K$ with $z$-transform $z$.
Clearly $T_0\seq T|_K$. By part~4 of proposition~\ref{t1.1}, 
$T_0=T|_K$.\qed

\bppsn\label{t1.4}
Let $S$ be a closed semiregular operator on $\cla$ such that $S|_K$
is regular with $z$-transform $z\in M(K)$. Then for any $a\in D(S)$,
there is a $c\in M(K)$ such that 
\bean
a &=& (I-{z}^*z)^{1/2}c,\\
Sa &=& z c.
\eean
\eppsn
\NI {\it Proof\/}: 
Take an $a\in D(S)$. Let
$\{e_\afa\}$ be an approximate identity for $K$.
For each $\afa$, one has $ae_\afa\in D(S)\cap K=D(S|_K)$.
Hence there is a $c_\afa \in K$ such that
\be \label{e1}
\brray{rcl}
ae_\afa &=& (I-{z}^*z)^{1/2}c_\afa,\\
S(ae_\afa) &=& z c_\afa.
\erray
\ee
From the above equations it follows that 
$c_\afa=(I-{z}^*z)^{1/2}ae_\afa+{z}^*(Sae_\afa)=ce_\afa$, 
where  $c=(I-{z}^*z)^{1/2}a+{z}^*(Sa)$.
Now using the fact that $e_\afa$ is an approximate identity, 
we get
\bean
ak &=& (I-{z}^*z)^{1/2}ck,\\
(Sa)k &=& z c k
\eean
for all $k\in K$, which proves the result.\qed

The above proposition together with the one that
follows will be the key ingredients in proving the regularity
of certain semiregular operators later.

\bppsn\label{t1.5}
For any $a\in D(S^*)$, there exists $c\in M(K)$ such that
\bean
a &=& (I-z{z}^*)^{1/2}c,\\
S^*a &=& {z}^* c.
\eean
\eppsn
\prf Similar to the proof of the previous proposition.\qed

Let us denote by $D$ the set $\{(I-{z}^*z)^{1/2}a+{z}^*(Sa):a\in D(S)\}$
and by $D_*$ the set $\{(I-z{z}^*)^{1/2}a+z(S^*a):a\in D(S^*)\}$. Observe 
that for $c\in D$ and $d\in D_*$, $zc$ and $z^*d$ are in $\cla$.

\blmma\label{t1.6}
Let $D$ be as above, and assume that $S=S^{**}$. Then
\begin{enumerate}
\item $D$ is a Hilbert $\cla$-module contained in $M(K)$,
\item $D=\Gamma(z):=(I-{z}^*z)^{-1/2}\cla \cap {z}^{-1}\cla
        \equiv\{c\in M(K): (I-{z}^*z)^{1/2}c\in \cla, zc\in \cla\}.$
\end{enumerate}
\elmma
\prf
Part~1 is straightforward. We will prove part~2 here.
Define an operator $\tilde{S}:(I-{z}^*z)^{1/2}\Gamma(z)\raro \cla$ by
\[
\tilde{S}((I-{z}^*z)^{1/2}c)=zc,\quad c\in\Gamma(z).
\]
By proposition~\ref{t1.4}, $D\seq\Gamma(z)$ and $S\seq\tilde{S}$.
Hence $\tilde{S}$ is densely defined. From the injectivity of 
$(I-{z}^*z)^{1/2}$ it follows that $\tilde{S}$ is well-defined.
It can easily be verified from the definition of $\tilde{S}$ that
it is closed. 

By proposition~\ref{t1.1}, $S^*|_K=(S|_K)^*$ and hence has 
$z$-transform ${z}^*$. From proposition~\ref{t1.5}, we conclude that
$D_*\seq\Gamma({z}^*)$.
Now, for $d\in D_*$ and 
$c\in \Gamma(z)$,
\bean
\lgl (I-z{z}^*)^{1/2}d, \tilde{S}((I-{z}^*z)^{1/2}c) \rgl
   &=& \lgl (I-z{z}^*)^{1/2}d,zc \rgl \\
 &=& \lgl {z}^*(I-z{z}^*)^{1/2}d,c \rgl \\
 &=& \lgl S^*((I-z{z}^*)^{1/2}d),(I-{z}^*z)^{1/2}c \rgl,
\eean
so that $D(S^*)\seq D((\tilde{S})^*)$. Therefore $S^*\seq (\tilde{S})^*$.
Thus $S\seq\tilde{S}\seq (\tilde{S})^{**}\seq S^{**}=S$.
This implies $D(S)=D(\tilde{S})\equiv (I-{z}^*z)^{1/2}\Gamma(z)$,
i.e.\ $\Gamma(z)\seq D$.
\qed

A similar statement about $D_*$ also holds; except that in that case
one need not assume $S^*=S^{***}$, it is automatic.
The above proposition tells us that if $S|_K$ is regular, 
even though $S$  may not be regular, it is uniquely determined
by a bounded adjointable operator on $K$, as long as $S$ is 
sufficiently nice (i.e.\  $S=S^{**}$).

\bppsn\label{t1.7}
Let $S$ be a closed semiregular operator on $\cla$ such that 
$S|_K$ is regular with $z$-transform $z$. 
Then one has the following inclusions:
\[
\brray{rlccrl}
{\rm i.} & z\cla\seq\overline{(I-zz^*)^{1/2}\cla}, &&&
          {\rm ii.} & z^*\cla\seq\overline{(I-z^*z)^{1/2}\cla}, \\
{\rm iii.} & \cla z\seq\overline{\cla (I-z^*z)^{1/2}}, &&&
          {\rm iv.} & \cla z^*\seq\overline{\cla (I-zz^*)^{1/2}},\\
{\rm v.} & z^*z\cla\seq\overline{(I-z^*z)\cla}, &&&
          {\rm vi.} & zz^*\cla\seq\overline{(I-zz^*)\cla},\\
{\rm vii.} & \cla\seq\overline{(I-z^*z)\cla}, &&&
          {\rm viii.} & \cla\seq\overline{(I-zz^*)\cla}.
\erray
\]
{\rm (}here overline indicates closure in the norm topology\/{\rm )}
\eppsn
\prf We will prove (i) here. Proof of (ii) is similar. All the
other inclusions follow from these two.

Take any $a=(I-z^*z)^{1/2}d\in D(S)$. Then 
$za=z(I-z^*z)^{1/2}d=(I-zz^*)^{1/2}zd\in (I-zz^*)^{1/2}\cla$.
Thus $zD(S)\seq (I-zz^*)^{1/2}\cla$. Since $D(S)$ is dense 
in $\cla$, we have the required inclusion. \qed

\bcrlre \label{t1.8}
With the notation as above, one has the following
\bean
 D &\seq & \overline{(I-z^*z)^{1/2}\cla},\\
 D_* &\seq &\overline{(I-zz^*)^{1/2}\cla}.
\eean
\ecrlre
\prf Any $d\in D$ is
of the form $(I-z^*z)^{1/2}a+z^*Sa$ for some $a\in D(S)$.
By part~(ii) of the previous proposition, 
$z^*Sa\in \overline{(I-z^*z)^{1/2}\cla}$. Hence we have the
first inclusion. Proof of the other one  is similar.\qed

\blmma \label{t1.9}
Let $S$ be as in proposition~\ref{t1.7}. If
$z\in M(\cla)$ then $S^{**}$ is regular.
\elmma

\prf  From corollary~\ref{t1.8} and the given condition,
it follows that $D\seq \cla$. Therefore $(I-z^*z)^{1/2}\cla$
contains $D(S)$ and is dense in $\cla$. So $z$ is indeed the
$z$-transform of some regular operator $T$ on $\cla$.
Clearly $S\seq T$, so that $T^*\seq S^*$. 
From corollary~\ref{t1.8} we also have $D_*\seq \cla$.
Therefore
$D(S^*)=(I-zz^*)^{1/2}D_* \seq (I-zz^*)^{1/2}\cla=D(T^*)$.
It follows then that $S^*=T^*$. Hence $S^{**}=T^{**}=T$.
Thus  $S^{**}$ is regular.\qed

\bppsn\label{t1.10}
Let $S$ and $z$ be as in the previous proposition. If 
$z^*z\in M(\cla)$ then $S^{**}$ is regular.
\eppsn

\prf Let us first show that $zz^*$ is also in $M(\cla)$.
Take any $a$ and $b$ in $D(S^*)$. There are elements
$c$, $d$ in $D_*$ such that $a=(I-zz^*)^{1/2}c$
and $b=(I-zz^*)^{1/2}d$. For any integer $n\geq 1$, we have
$a^*(zz^*)^n b=c^*(I-zz^*)^{1/2}z(z^*z)^{n-1}z^*(I-zz^*)^{1/2}d
              =(z^*c)^*(I-z^*z)^{1/2}(z^*z)^{n-1}(I-z^*z)^{1/2}z^*d \in\cla$.
Since $D(S^*)$ is norm dense in $\cla$, one has
$a^*(zz^*)^n b\in \cla$ for all $a,b\in\cla$. Which means in particular that
$zz^*$ and $(zz^*)^2$ both are in $QM(\cla)$, the space of 
quasi-multipliers of $\cla$. By proposition~5.3 in \cite{wo2},
$zz^*\in LM(\cla)$, and since $zz^*$ is positive, it is actually
in $M(\cla)$.

Now from parts~(i) and (iii) of proposition~\ref{t1.7} and the forgoing 
proposition, it follows that $S^{**}$ is regular.\qed

\newsection{Regularity}
We are now ready for the main results in this paper.
Let $\pi$ be the canonical projection of $M(K)$ onto
$M(K)/K$. Restriction of $\pi$ to $\cla$ gives the canonical 
projection of $\cla$ onto $\cla/K$.

\bthm \label{t1.20}
Let  $S$ be a closed semiregular operator on $\cla$ such
that  its restriction to $K$ is regular. If
\be\label{e2}
\Bigl(Z(\cla/K)\cap \pi(D(S))\Bigr)\cla/K \;\;\mbox{\it is total in } \cla/K,
\ee
where $Z(\cla/K)$ is the centre of $\cla/K$, then $S^{**}$ is regular.
\ethm

\prf  Let $z$ be the $z$-transform of 
$S|_K$, and let $\{e_{\afa}\}_\afa$ be an approximate identity
in $\cla$. By part~(iii) of proposition~\ref{t1.7}, there 
exist elements $f_\afa$ in $\cla$ such that 
\(
\lim_\afa \|e_\afa z - f_\afa(I-z^*z)^{1/2}\|=0.
\)
This implies that 
\[
\lim_\afa \|z^*e_\afa^2 z -(I-z^*z)^{1/2}{f_\afa}^* f_\afa(I-z^*z)^{1/2}\|=0,
\]
which, in turn, implies that
\[
\lim_\afa \|z^*zd -(I-z^*z)^{1/2}{f_\afa}^* f_\afa(I-z^*z)^{1/2}d\|=0
\]
for all $d\in D$. It follows then that
\[
\lim_\afa \|(I-z^*z)^{1/2}d -(I-z^*z)(I+{f_\afa}^* f_\afa)(I-z^*z)^{1/2}d\|=0
\]
for all $d\in D$, i.e.\ 
\[
\lim_\afa \|a -(I-z^*z)(I+{f_\afa}^* f_\afa)a\|=0 \quad \forall a\in D(S).
\]
Applying $\pi$ now, we get
\[
\lim_\afa \|\pi(a) -(I-\pi(z)^*\pi(z))(I+{\pi(f_\afa)}^* 
       \pi(f_\afa))\pi(a)\|=0 \quad \forall a\in D(S).
\]
Now choose an $a\in D(S)$ such that $\pi(a)\in Z(\cla/K)$, 
then $I+\pi(f_\afa)^*\pi(f_\afa)$ will commute with $\pi(a)$.
Therefore using  the facts that
$\|(I+{f_\afa}^*f_\afa)^{-1}\|\leq 1$ and 
$(I+\pi(f_\afa)^*\pi(f_\afa))^{-1}$ also commutes with $\pi(a)$,
we get
\be \label{e1.1}
\lim_\afa \|(I+{\pi(f_\afa)}^* \pi(f_\afa))^{-1}\pi(a)
 -(I-\pi(z)^*\pi(z))\pi(a)\|=0
\ee
for all $\pi(a)\in Z(\cla/K)\cap \pi(D(S))$. From condition~(\ref{e2}),
it follows that (\ref{e1.1}) holds for all $\pi(a)\in \cla/K$.
That is, for any $a\in\cla$, $\pi(z)^*\pi(z)\pi(a)\in \cla/K$.
Hence there is a $b\in\cla$ and a $k\in K$ such that 
$z^*za=b+k$, which implies that $z^*za\in\cla$. Thus $z^*z\in M(\cla)$.
From proposition~\ref{t1.10}, we conclude that $S^{**}$ is regular.\qed

The following two corollaries are now immediate.

\bcrlre \label{t1.11}
Let $S$ be a closed semiregular operator on $\cla$ such that
its restriction to $K$ is regular. If $\cla/K$
is abelian, then $S^{**}$ is regular.
\ecrlre
\prf In this case, $Z(\cla/K)\cap \pi(D(S))=\pi(D(S))$.
Therefore condition~(\ref{e2}) holds.\qed

\bcrlre \label{t1.12}
Let $S$ be as in the earlier theorem. If $\cla/K$ is unital,
then $S^{**}$ is regular.
\ecrlre
\prf Since $\pi(D(S))$ is a dense right ideal in $\pi(\cla)=\cla/K$
which is unital, we have $\pi(D(S))=\cla/K$. Therefore 
$I\in Z(\cla/K)\cap \pi(D(S))$. So (\ref{e2}) is satisfied.\qed

\brmrk \label{t1.13}
{\rm
 We will primarily be interested in the case $K=\clb_0(\clh)$.
By proposition~5.1 of \cite{pa2}, the condition that the restriction
of $S$ to $K$ is regular is automatic in this case. }
\ermrk

It is now natural to ask what happens in the general case, i.e.\ 
when $\cla/K$ is neither unital nor abelian. We will give a 
counterexample to illustrate that the result may fail to hold
in general. Before going to the example, let us observe that
if $S$ is a semiregular operator on $\cla$, then the prescription
\bean
D(\pi(S)) &:= & \pi(D(S)), \\
\pi(S)\pi(a) &:=&\pi(Sa),\quad a\in D(S),
\eean
defines a semiregular operator on $\pi(\cla)$. The example below,
which appears in \cite{hil} as an example of a nonregular operator,
will in fact show that even if $S|_K$ and $\pi(S)$ both are 
regular, $S$ may fail to be so.

Let us first define an operator on the Hilbert $C^*$-module
$E=C[0,1]\otimes\clh$, where $\clh=L_2(0,1)$. Let 
 $\bta$ be the following function on the interval $[0,1]$:
\[
\bta(\pi)=\cases{1 & if $\pi=0$,\cr
                 \exp(i/\pi) & if $0<\pi\leq 1$.}
\]
Let 
\[
D_\pi = \{f\in L_2(0,1):\; f \mbox{ absolutely continuous, }
           f'\in L_2(0,1), f(0)=\bta(\pi)f(1)\},
\]
For $f\in E$, denote by $f_\pi$ the function $f(\pi,\cdot)$ in
$\clh$. Let $T$ be the semiregular operator on $E$ defined as follows:
\[
D(T)=\{f\in E : f_\pi\in D_\pi\; \forall\, \pi, \pi\mapsto (f_\pi)' 
                        \mbox{ is continuous}\}
\]
\[
(Tf)_\pi := i(f_\pi)'.
\]
It has been shown by Hilsum~(\cite{hil}) that this is a self-adjoint
nonregular operator. Also, from proposition~2.9 in \cite{hil}, it
follows that the restriction of $T$ to the submodule 
$F=C_0(0,1]\otimes\clh$ is a self-adjoint regular operator.

Notice two things now. $\cla=C[0,1]\otimes \clb_0(\clh)$ is the 
$C^*$-algebra of `compact' operators on $E$, and 
$K=C_0(0,1]\otimes \clb_0(\clh)$ is the corresponding $C^*$-algebra
for $F$. $K$ can easily be seen to be an essential ideal in $\cla$,
and $\cla/K \cong \clb_0(\clh)$. Let $\phi_1$ be the map introduced
in section~3 of \cite{pa2} for the Hilbert module $E$. Define $S$
to be the operator $\overline{\phi_1(T)}$ on $\cla$.
Using lemmas 3.1, 3.2 and~3.5 in \cite{pa2}, one can prove that for
any semiregular operator $\bbt$ on $E$, 
$\overline{\phi_1(\bbt^*)}={\phi_1(\bbt)}^*$. Since in our case
$T$ is self-adjoint, it follows that $S$ is self-adjoint too.
Nonregularity of $S$ is also clear by the discussion at the end
of section~3 in \cite{pa2}. Restriction of $S$ to $K$ is the
$\phi_1$-image of the restriction of $T$ to $F$. Therefore
$S|_K$ is regular. Since $\cla/K\cong \clb_0(\clh)$, the projection
of $S$ on $\cla/K$ is also regular by proposition~5.1 in \cite{pa2}.

\brmrk
{\rm
If we write $z$ for the $z$-transform of the restriction of $S$
to $K$, then the above example tells us that the inclusions in 
proposition~\ref{t1.7} are not enough to guarantee that $z\in M(\cla)$,
as in that case $S$ would have been regular.}
\ermrk

\newsection{Examples}
We will restrict ourselves to two examples in this section 
that occur naturally in the study of quantum groups.
The first one is the $C^*$-algebra corresponding to the quantum complex
plane and the other one is the crossed product algebra
$C_0(q^{\IZ}\cup \{0\})\ltimes_\afa \IZ$, where $q$ is a fixed real number in the
interval (0,1), $q^\IZ$ stands for the set $\{q^k:k\in\IZ\}$, and the
action $\afa$ of $\IZ$ on $C_0(q^\IZ)$ is given by
\bean
\afa_k f(q^r) &=& f(q^{r-k}),\quad r,k\in\IZ,\\
\afa_k f(0) &=& f(0).
\eean

Let us start with the quantum complex plane. Let $\clh=L_2(\IZ)$,
with canonical orthonormal basis $\{e_n\}_n$. Let $\ell^*$ and $q^N$
denote the following operators:
\[
\brray{rcll}
\ell^* e_k & = & e_{k+1}, & k\in\IZ,\\
q^N e_k & = & q^k e_k, & k\in\IZ.
\erray
\]
Let $D$ denote the linear span of
 $\{(\ell^*)^k f_k(q^N): k\in\IZ, f_k\in C_0(q^\IZ\cup\{0\}),
                                  f_k(0)=0 \mbox{ for } k\neq 0\}$.
The $C^*$-algebra of `continuous vanishing-at-infinity functions'
on the quantum plane, which we denote by $C_0(\IC_q)$, is the 
norm closure of $D$. The quantum complex plane can be looked upon
as the  homogeneous space $E_q(2)/S^1$ ($S^1$ being the one dimensional 
torus) for the quantum $E(2)$ group (\cite{pa1},\cite{wo1}).
$C_0(\IC_q)$ was introduced in a slightly different form
in \cite{wo1} (For a proof of the fact that the $C^*$-algebra described
above is isomorphic to the one in \cite{wo1}, see \cite{pa1}).

\blmma\label{t1.14}
$C_0(\IC_q)/\clb_0(\clh)\cong\IC$.
\elmma

\prf It is easy to see that $C_0(\IC_q)$ acts irreducibly on $\clh$
and contains the compact operator $|e_0\rgl\lgl e_0|= I_{\{1\}}(q^N)$.
Therefore $\clb_0(\clh)\seq C_0(\IC_q)$.

Define a map $\phi:C_0(\IC_q)\raro \IC$ by the prescription
\[
\phi\biggl(\sum_k (\ell^*)^k f_k(q^N)\biggr)=f_0(0),
          \quad \sum_k (\ell^*)^k f_k(q^N)\in D.
\]
It extends to a complex homomorphism of $C_0(\IC_q)$. 
It is easy to see that $\ker \phi$ is the closure of
$\{(\ell^*)^k f_k(q^N): k\in\IZ, f_k\in C_0(q^\IZ\cup\{0\}),
                                  f_k(0)=0 \mbox{ for all } k\}$,
i.e.\ is isomorphic to $C_0(\IZ)\ltimes \IZ$, which in turn is isomorphic
to $\clb_0(\clh)$.\qed

We can now apply corollary~\ref{t1.12} to conclude that for 
any closed semiregular operator $S$ on $C_0(\IC_q)$, $S^{**}$ is regular. 
Indeed, since the restriction of
$S$ to $\clb_0(\clh)$ is regular, by proposition~\ref{t1.2}, 
$S^{**}$ is an operator satisfying the assumptions of corollary~\ref{t1.12}.
\vspace{1.5ex}

Our second example, the crossed product algebra 
$\cla=C_0(q^\IZ \cup \{0\})\ltimes\IZ$,
is actually very similar to the previous one. Its relevance in quantum groups
stems from the fact that for any infinite dimensional irreducible representation
$\pi$ of the $C^*$-algebra $C_0(E_q(2))$ corresponding to the quantum $E(2)$ 
group, $\pi(C_0(E_q(2)))$ is isomorphic to $\cla$. From the definition of a
crossed product algebra, it can be shown quite easily that
$\cla$ is the norm closure of the linear span of 
 $\{(\ell^*)^k f_k(q^N): k\in\IZ, f_k\in C_0(q^\IZ\cup\{0\})\}$.
One then shows that $\cla/\clb_0(\clh)\cong C(S^1)$. The proof
is similar to the proof of lemma~\ref{t1.14}, except that the map
$\phi$ in this case maps $\cla$ onto $C(S^1)$ and is defined by
$\phi(\sum_k (\ell^*)^k f_k(q^N))=\sum_k f_k(0)\zeta^k$,
where $\zeta$ stands for the function $z\mapsto z$ on $S^1$.




\begin{thebibliography}{99}
\bibitem{dix} Dixmier, J. : {\sl $C^*$-Algebras}, North-Holland, 1977.
\bibitem{hil} Hilsum, M. : Fonctorialit\'{e} en K-th\'{e}ory bivariante
    pour les vari\'{e}t\'{e}s lipschitziennes,
     {\em $K$-Theory}, 3(1989), 401--440.
\bibitem{la} Lance, E. C. : {\sl Hilbert $C^*$-modules - A Toolkit for Operator
    Algebraists}, Cambridge University Press, 1995.
\bibitem{pa1} Pal, A. : {\sl On Some Quantum Groups and Their
    Representations}, Ph.\ D. Thesis, Indian Statistical Institute, 1995.
\bibitem{pa2} Pal, A. : Regular operators on Hilbert $C^*$-modules,
   {\em Preprint}, 1997, (to appear in the {\em Journal of Operator Theory}).
\bibitem{pe} Pedersen, G.K. : {\sl $C^*$-algebras and Their 
      Automorphism Groups}, Academic Press, 1979.
\bibitem{wo1} Woronowicz, S.L. : Unbounded Elements Affiliated With
  $C^*$-algebras and Noncompact Quantum Groups,
    {\em Comm. Math. Phys.}, 136(1991), 399--432.
\bibitem{wo2} Woronowicz, S.L. :   $C^*$-algebras generated
  by unboounded elements,
    {\em Rev. Math. Phys.}, 7(1995), No.\ 3, 481--521.

\end{thebibliography}
\end{document}